\documentclass[a4paper, 12pt]{article}

\usepackage{color}
\usepackage[dvips]{graphicx}
\usepackage{amscd}
\usepackage{amsmath}
\usepackage{ascmac}
\usepackage{cases}
\usepackage{amssymb}
\usepackage{amsfonts}
\usepackage{amsthm}
\usepackage[all]{xy}
\usepackage{enumerate}
\usepackage{wrapfig}
\usepackage{bm}
\usepackage{color}
\usepackage{comment}
\usepackage{authblk}

\theoremstyle{definition}
\newtheorem{Thm}{{\bf Theorem}}[section]

\newtheorem{Lem}[Thm]{{\bf Lemma}}
\newtheorem{Prop}[Thm]{{\bf Proposition}}
\newtheorem{Cor}[Thm]{{\bf Corollary}}

\newtheorem{Def}[Thm]{{\bf Definition}}

\newcounter{Exami}

\newcommand{\ilim}[1][]{\mathop{\varinjlim}\limits_{#1}}

\newcommand{\bnu}{{\boldsymbol \nu}}

\newcommand{\bP}{{\mathbb P}}      
\newcommand{\bZ}{{\mathbb Z}}

\newcommand{\cO}{{\mathcal O}}

\newcommand{\cM}{{\mathcal M}}
\newcommand{\cN}{{\mathcal N}}

\newcommand{\cK}{{\mathcal K}}

\newcommand{\cF}{{\mathcal F}}
\newcommand{\cE}{{\mathcal E}}

\newcommand{\cL}{{\mathcal L}}

\newcommand\Ext{\mathop{\rm Ext}\nolimits}

\newcommand\Pic{\mathop{\rm Pic}\nolimits}
\newcommand\Spec{\mathop{\rm Spec}\nolimits}
\newcommand\Hilb{\mathop{\rm Hilb}\nolimits}

\newcommand{\res}{\mathop{\sf res}\nolimits}

\newcommand{\codim}{\mathop{\rm codim}\nolimits}
\newcommand{\supp}{\mathop{\rm supp}\nolimits}

\newcommand{\depth}{\mathop{\rm depth}\nolimits}

\newcommand\lra{\longrightarrow}
\newcommand\ra{\rightarrow}

\DeclareMathOperator{\Bl}{Bl}

\makeatletter
\newcommand{\doublewidetilde}[1]{{%
  \mathpalette\double@widetilde{#1}%
}}
\newcommand{\double@widetilde}[2]{%
  \sbox\z@{$\m@th#1\widetilde{#2}$}%
  \ht\z@=.9\ht\z@
  \widetilde{\box\z@}%
}
\makeatother

\newcommand{\Address}{{
  \bigskip
  \footnotesize
  
  Yuki Matsubara, \textsc{Department of Mathematics, Graduate School of Science, Kobe University, 1-1 Rokkodaicho, Nada-ku, Kobe, 657-8501, Japan}\par\nopagebreak
  \textit{E-mail address}, Yuki Matsubara: \texttt{ymatuba@math.kobe-u.ac.jp}
}}

\begin{document}

\title{On the Cohomology of the Moduli Space of Parabolic Connections}
\author{Yuki Matsubara}
\date{}
\maketitle
\begin{abstract}
We study the moduli space of logarithmic connections of rank $2$ on $\mathbb{P}^1 \setminus \{ t_1, \dots, t_5 \}$ with fixed spectral data.
The aim of this paper is to compute the cohomology of this space,
a computation that will be used to extend the results of the Geometric Langlands Correspondence due to D. Arinkin to the case where these types of connections have five simple poles on $\mathbb{P}^1$.
\end{abstract}

\section{Introduction}
In this paper, we study the moduli space of logarithmic connections of rank $2$ on $\mathbb{P}^1 \setminus \{ t_1, \dots, t_n \}$ with fixed spectral data.
These moduli spaces have been studied from various points of view.
For example, they occur as spaces of initial conditions for Garnier systems (\cite{IIS}).
In a recent paper \cite{S},
C.\ Simpson studied some of the topological structures of related moduli spaces in the context of problems such as the WKB theories, and the $P = W$ conjecture. 
Our interest in the subject of moduli spaces comes from its relation to the Geometric Langlands Correspondence.
In \cite{A}, D.\ Arinkin proved such correspondence in a special case, by using the geometry of the moduli space of such connections on $\mathbb{P}^1 \setminus \{ t_1, \dots, t_4 \}$.
If $n \geq 5$, this moduli space has not been studied in detail, for its dimension is $2(n-3)$, which is larger than $4$.
In this work, by using canonical coordinates introduced by apparent singularities, we are able to reduce the problem to that of the geometry of surfaces (see \S \ref{Apparent singularities}).

\subsection*{The logarithmic connections.}
Fix points $t_1, \dots, t_n \in \mathbb{P}^1 (t_i \neq t_j)$, and set $D = t_1 + \cdots + t_n$.
We consider pairs $(E, \nabla)$, where $E$ is a rank $2$ vector bundle on $\mathbb{P}^1$ and $\nabla : E \ra E \otimes \Omega^1_{\mathbb{P}^1}(D)$ is a connection having simple poles supported by $D$.
At each pole, we have two residual eigenvalues $\{ \nu_i^+, \nu_i^- \}$ of $\nabla$ for each $i = 1, \dots, n$; they satisfy the Fuchs relation $d + \sum_i(\nu_i^+ + \nu_i^-)= 0$, where $d := \deg(E)$.
Moreover, we can naturally introduce parabolic structures $\mbox{\boldmath $l$} = \{l_i \}_{1 \leq i \leq n}$ such that $l_i$ is a one-dimensional subspace of $E_{t_i}$ which corresponds to an eigenspace of the residue of $\nabla$ at $t_i$ with the eigenvalue $\nu_i^+$.
Note that, when $\nu_i^+ \neq \nu_i^-$, the parabolic structure $\mbox{\boldmath $l$}$ is determined by the connection $(E, \nabla)$.
Fixing spectral data $\bnu = (\nu_i^{\pm})$ with integral sum $-d$, by introducing the weight $\mbox{\boldmath $w$}$ for stability, one can construct the moduli space $M^{\mbox{\boldmath $w$}}(\mbox{\boldmath $t$}, \bnu, d)$ of $\mbox{\boldmath $w$}$-stable $\bnu$-parabolic connections $(E, \nabla, \mbox{\boldmath $l$})$ of degree $d$ using Geometric Invariant Theory, and the moduli space $M^{\mbox{\boldmath $w$}}(\mbox{\boldmath $t$}, \bnu, d)$ turns out to be a smooth irreducible quasi-projective variety of dimension $2(n-3)$ (see \cite{IIS} for details). 

We note that, when $\sum_{i = 1}^n \nu_i^{\epsilon_i} \not\in \mathbb{Z}$, for any choice $(\epsilon_i) \in \{ +, - \}^n$, every parabolic connection $(E, \nabla, \mbox{\boldmath $l$})$ is irreducible, and thus stable for any weight $\mbox{\boldmath $w$}$; the moduli space $M^{\mbox{\boldmath $w$}}(\mbox{\boldmath $t$}, \bnu, d)$ does not depend on the choice of weights $\mbox{\boldmath $w$}$ in that case.

These moduli spaces occur as spaces of initial conditions for Garnier systems, the case $n = 4$ corresponding to the Painlev\'e VI equation.
Such differential equations are nothing but isomonodromic deformations for linear connections.
By suitable transformations, we may normalize $\bnu$ as 
\begin{equation*}
\begin{cases}
\nu_i^{\pm}=\pm \nu_i& (i=1,\ldots,n-1)\\
\nu_n^+=-d - \nu_n  \\
\nu_n^-= \nu_n,
\end{cases}
\end{equation*}
for some $(\nu_1,\ldots,\nu_n) \in \mathbb{C}^{n}$.
Denote by $\cM(d)$ the moduli stack of  $\bnu$-$\mathfrak{sl}_2$-parabolic connections of degree $d$ and by $M(d)$ its coarse moduli space.
By the above normalization, we have a natural isomorphism $M(d) \simeq M^{\mbox{\boldmath $w$}}(\mbox{\boldmath $t$}, \bnu, d)$ (see \cite{IIS}).
Moreover, $M(d)$ has a natural compactification $\overline{M(d)}$, which is the moduli space of $\lambda$-$\bnu$-parabolic connections $(E, \nabla_{\lambda}, \lambda \in \mathbb{C})$ over $\mathbb{P}^1$.
(Note that the moduli space $M(d)$ is nothing but the moduli space of $(\nu_1,\ldots,\nu_n)$-bundles on $\mathbb{P}^1$ treated in \cite{AL} and \cite{Obl}, and $\overline{M(d)}$ is the moduli space of $\epsilon$-bundles on $\mathbb{P}^1$ in \cite{A}).

We should mention that P.\ Boalch has a number of related works concerned with the case of meromorphic connections with irregular singularities.
We refer to \cite{B}, for example.

\subsection*{Main Results.}

\begin{Thm}\label{coh of M}

\textit{
Let  $\cM(d)$ be the moduli stack of  $\bnu$-$\mathfrak{sl}_2$-parabolic connections of degree $d$.
Then we have\\
$$H^i(\mathcal{M}(d), \mathcal{O}_{\mathcal{M}(d)}) =  \begin{cases}
                                                                                     \mathbb{C}, & i = 0, \\
                                                                                     0, & i > 0.
                                                                              \end{cases} $$}
\end{Thm}

\section{Preliminaries}
\subsection{$\mathfrak{sl}_2$-connections.}
We introduce $\mathfrak{sl}_2$-connections.

Fix complex numbers $\nu_1, \dots, \nu_n \in \mathbb{C}$.
Suppose that $\nu_1 \cdots \nu_n \neq 0$ and
\begin{equation*}\label{generic condition 2}
\sum^n_{i=1} \epsilon_i \nu_i \notin \bZ
\end{equation*}
for any $(\epsilon_i) \in \{+, - \}^n$.

\begin{Def}
A \textit{$\bnu$-$\mathfrak{sl}_2$-parabolic connection on $\bP^1$} is a triplet $(E, \nabla,\varphi)$
such that 
\begin{enumerate}
 \item[(1)] $E$ is a rank $2$ vector bundle of degree $d$ on $\bP^1$,
 \item[(2)] $\nabla\colon E \ra E\otimes \Omega^1_{\bP^1}(D)$ is a connection, where $D := t_1 + \cdots +t_n$,
 \item[(3)] $\varphi\colon \bigwedge^2 E \simeq \cO_{\bP^1}(d)$ is a horizontal isomorphism, 
 \item[(4)] the residue $\res_{t_i}(\nabla)$ of the connection $\nabla$ at $t_i$ has eigenvalues $\nu_i^{\pm}$ for each $i$ ($1\le i \le n$).
\end{enumerate}
We call $\bnu = (\nu_i^{\pm})_{1 \leq i \leq n}$ {\it local exponents}.
\end{Def}
There exists a one dimensional subspace $l_i \subset E_{t_i}$ on which $\res_{t_i}(\nabla)$ acts as multiplication by $\nu_i^+$.
For generic $\bnu$, the parabolic direction $l_i$ is nothing but the eigenspace for $\res_{t_i}(\nabla)$ with respect to $\nu_i^+$ so that the parabolic data $\mbox{\boldmath $l$} = \{ l_i \}$ is uniquely determined by the connection $(E, \nabla, \varphi)$ itself.

In this paper, it is enough to consider the case where $d = -1$.
By suitable transformations, we may put 
\begin{equation*}
 \nu^{\pm}_i := \pm \nu_i \ \ (i=1,\ldots, n-1 ),\ \nu^+_n:=1-\nu_n,\ \nu^-_n := \nu_n \rlap{.}
\end{equation*}

Denote by $\cM(d)$ the moduli stack of $\bnu$-$\mathfrak{sl}_2$-parabolic connections on $\bP^1$,
and by $M(d)$ its coarse moduli space.
This moduli space is a smooth, irreducible quasi-projective algebraic variety of dimension $2(n-3)$ (\cite[Theorem 2.1]{IIS}).
Recall that $\cM(d)$ has a natural compactification $\overline{\cM(d)}$ which is the moduli stack of $\lambda$-$\bnu$-parabolic connections $(E, \nabla_{\lambda}, \varphi, \lambda \in \mathbb{C})$ over $\mathbb{P}^1$.
(Note that, in \cite{A}, $\lambda$-$\bnu$-parabolic connections are called $\epsilon$-bundles.)
Then, under the condition that $(E, \nabla, \varphi)$ is irreducible, Arinkin showed that the moduli stack $\overline{\cM(d)}$ is a complete smooth Deligne-Mumford stack \cite[Theorem 1]{A}.
Moreover, he also showed that the $\lambda = 0$ locus $\cM(d)_H \subset \overline{\cM(d)}$, which is the moduli stack of parabolic Higgs bundles, is also a smooth algebraic stack.
On the other hand, as remarked in the proof of \cite[Proposition 7]{A}, the coarse moduli space $\overline{M(d)}$ corresponding to $\overline{\cM(d)}$ is not smooth: it has quotient singularities.
As for the possible smooth compactification by $\phi$-parabolic-connections, one may refer to \cite{IIS}.

\subsection{Lower and upper modifications.}
In this subsection, following \cite[\S 2]{Obl}, we describe the lower and upper modifications.
Let $E$ be an algebraic vector bundle on $\bP^1$ of rank $2$ and of degree $d$.
Fix a point $t \in \bP^1$.
Let $l \subset E_t$ be a one-dimensional subspace. 
\begin{Def}
We call  
\begin{equation*}
(t, l)^{\text{low}}(E) := \{ s \in E \mid s(t) \in l \}, \quad (t, l)^{\text{up}}(E) := (t, l)^{\text{low}}(E) \otimes \cO_{\mathbb{P}^1}(t)
\end{equation*}
\textit{the lower and upper modifications} of $E$, respectively.
\end{Def}
The lower and upper modifications provide the exact sequences
\begin{equation*}
0\lra (t, l)^{\text{low}}(E) \lra E \lra E_{t}/l \lra 0,
\end{equation*}
\begin{equation*}
0\lra E \lra (t, l)^{\text{up}}(E) \lra l \otimes \cO_{\mathbb{P}^1}(t) \lra 0,
\end{equation*}
respectively.
In other words, we change our bundle by rescaling the basis of sections in the neighborhood of a point $t$ as follows:
given a local decomposition $V=l\oplus l'$ of $E\simeq V \otimes \cO$,
we take the local basis $\{ s_1(z), s_2(z) \}$ with $l\otimes \cO \simeq \langle s_1(z) \rangle$ and $l'\otimes \cO \simeq \langle s_2(z) \rangle$.
Then the basis of the lower modification $(t, l)^{\text{low}}$ of the bundle is generated by the sections $\{ s_1(z),(z-x) s_{2}(z)\}$,
and the upper one $(t, l)^{\text{up}}$ is given by $\{ (z-t)^{-1}s_1(z), s_{2}(z)\}$.
Consequently, in the punctured neighborhood, we may represent the actions of the modifications by the following gluing matrices
\begin{equation*}
(t, l)^{\text{low}}=\left(
\begin{array}{ll}
1 & 0 \\
 0 & (z-t)
\end{array}
\right),\quad
(t, l)^{\text{up}}=\left(
\begin{array}{ll}
(z-t)^{-1}  & 0 \\
 0 & 1
\end{array}
\right).
\end{equation*}

For the parabolic bundle $(E, \mbox{\boldmath $l$})$, we recall the geometrical properties of these modifications.
Denote by $\mathbb{P}(E, \mbox{\boldmath $l$})$ the projectivization of the parabolic bundle $(E, \mbox{\boldmath $l$})$.
It consists of the projective bundle $\mathbb{P}E$ together with a parabolic point $l_i$ in the fiber $F$ of each $t_i$.
In this situation, the lower and the upper modifications of $E$ are birational transformations of the total space $\text{tot}(\mathbb{P}E)$: these are the blowing-ups of the point $l_i \in \mathbb{P}E$ followed by the contraction of the total transform $\widetilde{F}$ of the fiber $F$.
The point resulting from this contraction gives the new parabolic direction $l'_i$.
We recall their properties in the following proposition:

\begin{Prop}
\textit{
Let $(E, \mbox{\boldmath $l$})$ be a parabolic bundle over $(\mathbb{P}^1, \mbox{\boldmath $t$} = \{ t_i \})$. Then the parabolic bundle $(E', \mbox{\boldmath $l$}') = (t_i,l_i)^{\text{low}}(E)$ satisfies the following properties:}
\begin{enumerate}
  \item[(1)] \textit{$\det (E', \mbox{\boldmath $l$}') = \det (E, \mbox{\boldmath $l$}) \otimes \cO_{\mathbb{P}^1}(-t_i)$.}
  \item[(2)] \textit{If $L \subset E$ is a line subbundle passing through $l_i$, its image by $(t_i,l_i)^{\text{low}}$ is a subbundle $L' \simeq L$ of $(t_i,l_i)^{\text{low}}(E)$ not passing through $l'_i$.}
  \item[(3)] \textit{If $L \subset E$ is a line subbundle not passing through $l_i$, its image by $(t_i,l_i)^{\text{low}}$ is a subbundle $L' \simeq L \otimes \cO_{\mathbb{P}^1}(-t_i)$ of $(t_i,l_i)^{\text{low}}(E)$ passing through $l'_i$.}
 \end{enumerate}
 \textit{For the upper modification, the parabolic bundle $(E'', \mbox{\boldmath $l$}'') = (t_i,l_i)^{\text{up}}(E)$ satisfies:}
  \begin{enumerate}
  \item[(4)] \textit{$\det (E'', \mbox{\boldmath $l$}'') = \det (E, \mbox{\boldmath $l$}) \otimes \cO_{\mathbb{P}^1}(t_i)$.}
  \item[(5)] \textit{If $L \subset E$ is a line subbundle passing through $l_i$, its image by $(t_i,l_i)^{\text{up}}$ is a subbundle $L' \simeq L \otimes \cO_{\mathbb{P}^1}(t_i)$ of $(t_i,l_i)^{\text{up}}(E)$ not passing through $l''_i$.}
  \item[(6)] \textit{If $L \subset E$ is a line subbundle not passing through $l_i$, its image by $(t_i,l_i)^{\text{up}}$ is a subbundle $L' \simeq L$ of $(t_i,l_i)^{\text{up}}(E)$ passing through $l''_i$.}
 \end{enumerate}
\end{Prop}

For a $\bnu$-$\mathfrak{sl}_2$-parabolic connection $(E, \nabla, \varphi)$, the lower modification of $E$ gives the new connection $\nabla'$ which is deduced from the action of $\nabla$ on the subsheaf $(t_i,l_i)^{\text{low}}(E) \subset E$, and,  over $t_i$, local exponents are changed by
$$ (\nu_i^+, \nu_i^-)' = (\nu_i^- + 1, \nu_i^+) \ \ \text{(and other} \ \nu_j^{\pm}\  \text{are left unchanged for} \ j \neq i). $$
The lower modufication gives us a morphism of moduli spaces $M(d) \ra M(d-1)$.
The upper modification defines the inverse map, and therefore, we have $M(d) \simeq M(d-1)$.

\subsection{Hirzebruch surfaces and the blowing-ups.}\label{Hirzebruch blowing up}
To describe the moduli space $M(-1)$, we introduce some blowing-ups of the Hirzebruch surface $\mathbb{F}_{n-2}$.
Put $L :=\Omega^1_{\mathbb{P}^1}(D)$.
We consider  the surface $\mathbb{F}_{n-2}$ as the total space of $\mathbb{P}(\cO_{\mathbb{P}^1} \oplus L)$.
Denote by $s_{\infty}$ the section defined by $L$.
$\mathbb{F}_{n-2} \setminus s_{\infty}$ is naturally identified with the total space of  $L$.
In particular, the affine part of the fiber $F_i$ over $t_i$ has the natural chart
$\res_{t_i} \colon F_i \setminus s_{\infty} \xrightarrow{\sim} \mathbb{C}$ given by the residue of sections of $L$.
We define two points $\hat{\nu}_i^{\pm} \in F_i$ by $\res_{t_i}(\hat{\nu}_i^{\pm}) = \nu_i^{\pm}$.

Denote by $\widetilde{\mathbb{F}_{n-2}} := \Bl_{\hat{\nu}_i^{\pm}} \mathbb{F}_{n-2}$ the blowing-up of $\mathbb{F}_{n-2}$ at $\hat{\nu}_i^{\pm}$ for each $i = 1, \dots, n$, by $\widetilde{s}_{\infty}, \widetilde{F}_i$ the strict transforms,
and by $E^{\pm}_i$ the exceptional curves at $(t_i, \hat{\nu}_i^{\pm})$.
Set
$$ \cK'_n := \widetilde{\mathbb{F}_{n-2}} \setminus (\widetilde{s}_{\infty} \cup \widetilde{F}_1 \cup \cdots \cup \widetilde{F}_n).$$
We denote by $\cK_n$ the image of $\cK'_n$ under the projection $\cK'_n \rightarrow \mathbb{F}_{n-2} \setminus s_{\infty}$.

\subsection{Apparent singularities and the dual parameters.}\label{Apparent singularities}
Let $(E,\nabla,\varphi) \in M(-1)$.
We can define the \textit{apparent singularities of $(E,\nabla,\varphi) \in M(-1)$} as follows:
we fix a section $s \in H^0(\bP^1, E)$.
For the section $s$, we define the following composition
\begin{equation*}
\cO_{\bP^1} \xrightarrow{\ s\ } E \xrightarrow{\ \nabla\ } E \otimes L \lra (E/\cO_{\bP^1}) \otimes L.
\end{equation*}
The composition $\cO_{\bP^1}\ra (E/\cO_{\bP^1}) \otimes L$ is an $\cO_{\bP^1}$-morphism, which is injective.
Then we can define a subsheaf $F^0\subset E $ such that $\cO_{\bP^1} \ra (F^0/\cO_{\bP^1}) \otimes L$ is an isomorphism.
By the isomorphism $F^0/\cO_{\bP^1} \simeq L^{-1}$, we have $F^0 \simeq \cO_{\bP^1} \oplus  L^{-1}$.
Therefore, we have the following exact sequence
\begin{equation*}\label{ES of App for conn}
0 \lra \cO_{\bP^1} \oplus L^{-1} \lra E \lra T_A \lra 0,
\end{equation*}
where $T_A$ is a torsion sheaf.
By the Riemann-Roch theorem, the torsion sheaf $T_A$ is length $n-3$.

\begin{Def}
For $(E,\nabla , \varphi) \in M(-1)$ and a nonzero section $s \in H^0(\bP^1,E)$, 
we call the support of $T_A$ {\it the apparent singularities of a $\bnu$-$\mathfrak{sl}_2$-parabolic connection with a cyclic vector $(E,\nabla , \varphi,[s])$}.
\end{Def}

Now, we consider the following stratification of  $M(-1)$. 
By the irreducibility of $(E, \nabla, \varphi) \in M(-1)$, we have the following proposition.
\begin{Prop}\label{bundle type}
\textit{
For $(E ,\nabla ,\varphi) \in M(-1)$,
we have 
\begin{equation*}
E\simeq \cO(k) \oplus\cO(-k-1) \text{  where } 0\le k \le \left[ \frac{n-3}{2} \right].
\end{equation*}
}
\end{Prop}
Denote by $M(-1)^k$ the subvariety of $M(-1)$ where $E\simeq \cO(k) \oplus\cO(-k-1)$.
Then
\begin{equation*}
M(-1) = M(-1)^0\cup \cdots \cup M(-1)^{[ (n-3)/2]}.
\end{equation*}
Note that the stratum $M(-1)^0$ is a Zariski open dense set of $M(-1)$.

For $(E,\nabla , \varphi) \in M(-1)^0$, we define \textit{dual parameters} as follows:
put $U_0 := \mathbb{P}^1 \setminus \{ \infty \}, U_{\infty} := \mathbb{P}^1 \setminus \{ 0 \}$.
Let $z$ and $w$ be coordinates on $U_0$ and $U_{\infty}$, respectively.
Put
$$ \omega_z := \frac{dz}{\prod_{i = 1}^n(z - t_i)} \quad
\text{and} \quad
R_0 :=\left(
\begin{array}{clcl}
 1 & 0 \\
 0 & \frac{1}{z}
\end{array}
\right). $$
Since $E \simeq \cO_{\mathbb{P}^1} \oplus \cO_{\mathbb{P}^1}(-1)$, we can denote the connection $\nabla$ by 
\begin{equation*}
\nabla = 
\begin{cases}
d+ A_z^0 \otimes \omega_z    & \text{on } U_0 \\
d+ R_0^{-1}dR_0 +  R_0^{-1}(A_z^0 \otimes \omega_z )R_0  & \text{on } U_{\infty},
\end{cases}\quad \text{where }\ 
A_z^0 :=\left(
\begin{array}{clcl}
f_{11}^{(n-2)}(z) & f_{12}^{(n-1)}(z) \\
 f_{21}^{(n-3)}(z) & -f_{11}^{(n-2)}(z)
\end{array}
\right).
\end{equation*}
Note that the zeros of the polynomial $f_{21}^{(n-3)}(z)$ are the apparent singularities of $(E,\nabla , \varphi)$.
We denote by $\{ q_1,\ldots,q_{n-3} \}$ the apparent singularities.
We put $p_i := f_{11}^{(n-2)}(q_i) \in L_{q_i}$.
We call $\{ p_1 ,\ldots, p_{n-3} \}$ the \textit{dual parameters} of $(E,\nabla , \varphi) \in M(-1)^0$.

\section{Geometric description of $M(-1)^0$}

Let $\cK_n'$ be the Zariski open set of the blowing-up of the Hirzebruch surface of degree $n-2$ defined in subsection \ref{Hirzebruch blowing up},
and let $\cK_n$ be the contraction $\cK'_n \ra \cK_n$. 
Then we can define the map
\begin{equation}\label{M0 to Sym}
\begin{aligned}
 M(-1)^0 &\lra \mathrm{Sym}^{n-3}(\cK_n) \\
 (E,\nabla, \varphi) &\longmapsto \{(q_1,p_1),\ldots,(q_{n-3},p_{n-3}) \},
\end{aligned}
\end{equation}
which was constructed in \cite[\S 3]{Obl}.
We consider the composite of the Hilbert-Chow morphism and the blowing-up
\begin{equation*}
\Hilb^{n-3}(\cK'_n) \lra \mathrm{Sym}^{n-3}(\cK'_n) \lra \mathrm{Sym}^{n-3}(\cK_n) \rlap{,}
\end{equation*}
where $\cK_n' \ra \cK_n$ is the blowing-up defined in subsection  \ref{Hirzebruch blowing up}.
We have the following proposition.

\begin{Prop}[\cite{KS} Theorem 5.2]\label{injective M0}
\textit{
We can extend the map (\ref{M0 to Sym}) to
 $$ M(-1)^0 \longrightarrow \Hilb^{n-3}(\cK'_n), $$
 and this map is injective.
} 
\end{Prop}

Suppose $n = 5$.
We denote by $Z \subset \mathrm{Sym}^2(\cK'_5)$ the proper pre-image of $\{ (q_1, p_1), (q_1, -p_1) \} \subset \mathrm{Sym}^2(\cK'_5)$ under  the blowing-up $\mathrm{Sym}^2(\cK'_5) \ra \mathrm{Sym}^2(\cK_5)$, and by $\widetilde{Z} \subset \Hilb^2(\cK'_5)$ the proper pre-image of $Z$ under the Hilbert-Chow morphism $\Hilb^2(\cK'_5) \ra \mathrm{Sym}^2(\cK'_5)$.
Denote by
\begin{equation}
 \widetilde{\Hilb}^2(\cK'_5) \ra \Hilb^2(\cK'_5)
\end{equation}
the blowing-up along $\widetilde{Z}$, and by $\widehat{Z}$ the strict transform of $\widetilde{Z}$. 
We also denote by $(\cK'_5 \times \cK'_5)^{\sim}$ the blowing-up of $\cK'_5 \times \cK'_5$ along the ideal $(q_1 - q_2, p_1 - p_2)$, and by $(\cK'_5 \times \cK'_5)^{\approx}$ the blowing-up of $(\cK'_5 \times \cK'_5)^{\sim}$ along the ideal $(q_1 - q_2, p_1 + p_2)$.
Then $\Hilb^2(\cK'_5) = (\cK'_5 \times \cK'_5)^{\sim}/\mathfrak{S}_2$ and $\widetilde{\Hilb}^2(\cK'_5) = (\cK'_5 \times \cK'_5)^{\approx}/\mathfrak{S}_2$.

Now, using the above description, we define another important blowing-up of the Hirzebruch surface $\mathbb{F}_3$.
Fix $q_1 \in \mathbb{P}^1 \setminus \{ t_1, \cdots, t_5 \}$ and define the fiber $F_6$ over $q_1$.
We denote by $(\mathbb{F}_3)^{\approx}$ the blowing-up of $\widetilde{\mathbb{F}_3}$ at two points $\{(q_1, p_1), (q_1, -p_1) \}$ 
(when $p_1 = p_2 = 0$, it blows up twice at $(q_1, 0)$).
Set
$$ \cK'_{5, q_1} := (\mathbb{F}_3)^{\approx} \setminus (\widetilde{s}_{\infty} \cup \widetilde{F}_1\cup \cdots \cup \widetilde{F}_6), $$
where $\widetilde{F}_6$ is the strict transform of $F_6$.
We denote by $E^{\pm}_6$ the exceptional curves at $(q_1, \pm p_1)$, and denote by $\cK_{5, q_1}$ the image of $\cK'_{5, q_1}$ under the projection $\cK'_{5, q_1} \rightarrow \mathbb{F}_{3} \setminus s_{\infty}$.

\section{Geometric description of $\cK'_{5, q}$}

In this section, for the sake of simplicity, we write $\cK'_{5, q}$ for $\cK'_{5, q_1}$.

\begin{Prop}\label{second coh of K}
\textit{
Let $\cF$ be any quasi-coherent sheaf on $\cK'_{5,q}$.
Then $H^i(\cK'_{5, q}, \cF) = 0$ for $i \geq 2$.
}
\end{Prop}

\begin{proof}
Let $Q$ be a projective line doubled at the six points $\{t_1, \dots, t_5, q \}$.
We can define a natural projection $\cK'_{5, q} \rightarrow Q$.
Moreover, this map is an affine bundle, thus it is an affine morphism.
\end{proof}

Set $D_{q} := 2\widetilde{s}_{\infty} + \widetilde{F}_1 + \cdots + \widetilde{F}_6$.
Then
\begin{equation}\label{deg N 0}
(D_q, D_q) = (D_q, \widetilde{s}_{\infty}) = (D_q, \widetilde{F}_i) = 0.
\end{equation}

We also have $K := K_{(\mathbb{F}_3)^{\approx}} = -2\widetilde{s}_{\infty} - \sum_{i = 1}^5 \widetilde{F}_i + E_6^+ + E_6^-$.
By the Riemann-Roch theorem, we have
$$ \chi(\cO_{D_q}) = -\frac{D_q(D_q + K)}{2} =-1. $$
This implies the following statement.

\begin{Prop}\label{prop11}
\textit{
Let $\cE$ be a locally free sheaf on $D_q$ of rank $r$. Then
$$ \chi(\cE) = 2\deg(\cE|_{\widetilde{s}_{\infty}}) + \sum_{i = 1}^6 \deg(\cE|_{\widetilde{F}_i}) - r. $$
}
\end{Prop}
\begin{proof}
This follows from the Riemann-Roch theorem for an embedded curve (cf. \cite[Chapter 2, Theorem 3.1]{BHPV}).
\end{proof}

\begin{Lem}\label{lem8}
\textit{
Let $\cE$ be a nontrivial invertible sheaf on $D_q$ such that 
$\deg (\cE |_{\widetilde{s}_{\infty}}) = 0$, and either $\deg ( \cE |_{\widetilde{F}_i} ) = 0$ for all i, or one of the numbers is $\deg ( \cE |_{\widetilde{F}_i} ) = -1$,  another one is $1$, and the remaining three equal zero.
Then $H^i(D_q, \cE) = 0$ for $i \neq 1$, and $H^1(D_q, \cE) = \mathbb{C}$.
}
\end{Lem}

\begin{proof}
By Proposition \ref{prop11}, we have $\chi(\cE) = -1$. Therefore, it is enough to prove that $H^0(D_q, \cE) = 0$.

Assume the converse. Let $f \in H^0 (D_q, \cE)$, $f \neq 0$. 
Now $\chi (\cE) = \chi ( \mathcal{O}_{D_q})$, and $\cE \not\simeq \mathcal{O}_{D_q}$, so $f$ is zero on one of the irreducible components of $D_q$. 
We take $\tilde{F}_1$ to be this component.

We may assume that $\deg (\mathcal{E}|_{\widetilde{F}_i} ) \leq 0$ for $i \neq 1$.
The closed subscheme $D'_q := \widetilde{s}_{\infty} + \sum_{i \neq 1} \widetilde{F}_i \subset D_q$ is reduced and connected. 
Besides this, $\cE|_{D'_q}$ has nonpositive degree on any irreducible component of $D'_q$.
Therefore, either $f|_{D'_q} = 0$, or $f|_{D'_q}$ has no zero.

In the second case, $f|_C \neq 0$, where $C \subset D_q$ is any irreducible component.
Therefore, $f \in \ker(H^0(D_q, \cE) \rightarrow H^0(D'_q, \cE))$.
In other words, $f \in H^0(D_q, \cE \otimes \mathcal{I}_{D'_q})$, where $\mathcal{I}_{D'_q} := \{ \tilde{f} \in \mathcal{O}_{D_q} \ |\  \tilde{f}|_{D'_q} = 0 \}$ is the sheaf of ideals of $D'_q$. 

We have $\mathcal{I}_{D'_q} = \mathcal{O}_{(\mathbb{F}_3)^{\approx}}(-D'_q)/\mathcal{O}_{(\mathbb{F}_3)^{\approx}}(-D_q)$, and $\supp \mathcal{I}_{D'_q} = \widetilde{s}_{\infty} + \widetilde{F}_1$.
Hence, $\deg (\mathcal{I}_{D'_q}|_{\widetilde{F}_1}) |_{\widetilde{F}_1}= \deg (\mathcal{O}_{(\mathbb{F}_3)^{\approx}}(-D'_q)|_{\widetilde{F}_1}) = -1$.
Therefore, $\deg (\cE \otimes \mathcal{I}_{D'_q}) = \deg (\cE|_{\widetilde{F}_1}) -1 \leq 0$.
In the same way, $\deg (\cE \otimes \mathcal{I}_{D'_q})|_{\widetilde{s}_{\infty}} = \deg (\cE|_{\widetilde{s}_{\infty}}) -1 = -1$.
Since $\cE \otimes \mathcal{I}_{D'_q}$ is an invertible sheaf on the connected reduced scheme $\widetilde{s}_{\infty} + \widetilde{F}_1$,
this implies $f \in H^0(D_q, \cE \otimes \mathcal{I}_{D'_q}) = 0$.
\end{proof}

Set $\Pic^{0}(D_q) := \{ \cE \in \Pic(D_q) | \deg (\cE | _{\widetilde{s}_{\infty}}) = 0, \ \deg ( \cE | _{\widetilde{F}_i}) = 0 \ \text{for all} \ i \}$.

\begin{Prop}\label{prop12}
\textit{
$$\Pic^0(D_q) \simeq \mathbb{A}^2\rlap{.}$$
}
\end{Prop}
\begin{proof}
Set $D_q^{red} := \widetilde{s}_{\infty} + \sum_{i=1}^{6} \widetilde{F}_i \subset D_q $. 
Then $\Pic^0(D_q) = \ker (\Pic (D_q) \rightarrow \Pic ( D_q^{red} ) ) $. 
Set $\mathcal{O}' := \ker ( \mathcal{O}_{D_q}^{*} \rightarrow \mathcal{O}_{D_q^{red}}^{*})$.
Then the exact sequence $ 0 \rightarrow \mathcal{O}' \rightarrow \mathcal{O}_{D_q}^{*} \rightarrow \mathcal{O}_{D_q^{red}}^{*} \rightarrow 1 $ defines an isomorphism
$H^1(D_q, \mathcal{O}') \xrightarrow{\sim} \Pic^0(D_q)$.
However,  $\mathcal{O}'$ is a locally free $\mathcal{O}_{\widetilde{s}_{\infty}}$-module which satisfies $\deg (\mathcal{O}') = -(\widetilde{s}_{\infty} , D_q^{red}) = -3$.
Hence $\Pic^0(D_q)$ is a $2$-dimensional $\mathbb{C}$-space.
\end{proof}

\begin{Lem}\label{lem10}
\textit{
The sheaf $\cN_{D_q} := \cO_{(\mathbb{F}_3)^{\approx}}(D_q)|_{D_q}$ is not trivial.
}
\end{Lem}
\begin{proof}
Assume the converse.
Let $\sigma \in H^0(D_q, \cN_{D_q})$ be a global section of $\cN_{D_q}$ with no zeros.
Since $(\mathbb{F}_3)^{\approx}$ is a smooth rational projective variety, $H^1((\mathbb{F}_3)^{\approx}, \cO_{(\mathbb{F}_3)^{\approx}}) = 0$,
and therefore $\sigma \in H^0(D_q, \cN_{D_q}) = H^0((\mathbb{F}_3)^{\approx}, \cO_{(\mathbb{F}_3)^{\approx}}(D_q)/\cO_{(\mathbb{F}_3)^{\approx}})$
can be lifted to $s \in H^0 ((\mathbb{F}_3)^{\approx}, \cO_{(\mathbb{F}_3)^{\approx}}(D_q))$.
Then $(s)$ is an effective divisor equivalent to $D_q$, and $C :=\supp (s) \subset \cK'_{5, q}$.
Denote by $C'$ the image of $C$ under the blowing-down $\cK'_{5, q} \ra \cK_{5, q}$.
Then $C' \sim 2s_{\infty} + \sum_{i = 1}^6 F_i$.

Now let $f(x, y)$ be  a local equation for $C'$ on some local chart.
Then we can write $f(x, y) = y^2 + a_1(x)y + a_2(x)$, where $\deg a_i(x) = 3i$.
By definiton, $C'$ passes through $(t_i, \hat{\nu}_i^+)$ and $(q, p_1)$ with multiplicity $1$.
Since we put $\hat{\nu}_i^{\pm} = \Pi_{t_i \neq t_j} (t_i - t_j)\nu_i^{\pm}$, where
$$\nu^{\pm}_i := \pm \nu_i \ \ (i=1,\ldots, 4 ),\ \nu^+_5:=1-\nu_5,\ \nu^-_5 := \nu_5\rlap{,}$$
by Vieta's formula, $a_1(x)$ satisfies $a_1(t_i) = 0$ for $i = 1, \dots, 4$ and $a_1(q) = 0$.
This implies $a_1(x) \equiv 0$. However, then $0 =  (1-\nu_5) + \nu_5 = 1$, which is a contradiction.

\end{proof}

\begin{Prop}\label{vanish of N}
\textit{ For $k \neq 0$,
$H^i(D_q, (\cN_{D_q})^{\otimes k}) = 0$ if $i \neq 1$ and $H^1(D_q, (\cN_{D_q})^{\otimes k}) = \mathbb{C}$.}
\end{Prop}

\begin{proof}
By (\ref{deg N 0}), we have $\cN_{D_q} \in \Pic^0(D_q)$.
Lemma $\ref{lem10}$ and Proposition $\ref{prop12}$ imply $(\cN_{D_q})^{\otimes k} \not\simeq \cO_{D_q}$ for $k \neq 0$.
Lemma $\ref{lem8}$ completes the proof.
\end{proof}

\begin{Cor}\label{coh of O}
\textit{
$$ H^i(\cK'_{5, q}, \cO_{\cK'_{5, q}}) =   \begin{cases}
                                                                \mathbb{C}, & i = 0, \\
                                                                 H_{D_q}^2((\mathbb{F}_3)^{\approx}, \cO_{(\mathbb{F}_3)^{\approx}}), &i =1 \\
                                                                 0, & i > 1.
                                                      \end{cases} $$
}
\end{Cor}

\begin{proof}
By local cohomology theory, we have the long exact sequence
\begin{equation*}
 \begin{split}
   0 &\ra H_{D_q}^0((\mathbb{F}_3)^{\approx}, \cO_{(\mathbb{F}_3)^{\approx}}) \ra H^0((\mathbb{F}_3)^{\approx}, \cO_{(\mathbb{F}_3)^{\approx}}) \ra H^0(\cK'_{5, q}, \cO_{\cK'_{5, q}}) \\
      &\ra H_{D_q}^1((\mathbb{F}_3)^{\approx}, \cO_{(\mathbb{F}_3)^{\approx}}) \ra H^1((\mathbb{F}_3)^{\approx}, \cO_{(\mathbb{F}_3)^{\approx}}) \ra H^1(\cK'_{5, q}, \cO_{\cK'_{5, q}}) \\
      &\ra H_{D_q}^2((\mathbb{F}_3)^{\approx}, \cO_{(\mathbb{F}_3)^{\approx}}) \ra H^2((\mathbb{F}_3)^{\approx}, \cO_{(\mathbb{F}_3)^{\approx}}) \ra H^2(\cK'_{5, q}, \cO_{\cK'_{5, q}}) \\
      &\ra H_{D_q}^3((\mathbb{F}_3)^{\approx}, \cO_{(\mathbb{F}_3)^{\approx}}) \ra 0.
 \end{split}
\end{equation*} 
and $H_{D_q}^i((\mathbb{F}_3)^{\approx}, \cO_{(\mathbb{F}_3)^{\approx}}) \simeq \ilim[k]\Ext^i (\cO_{kD_q}, \cO_{(\mathbb{F}_3)^{\approx}}) \simeq \ilim[k]H^{i-1}((\mathbb{F}_3)^{\approx}, \cN_{kD_q})$.
The statement follows from proposition \ref{vanish of N} and the rationality of $(\mathbb{F}_3)^{\approx}$.
\end{proof}

\subsection*{Special case: $q_1 \in \{ t_1, \dots, t_5 \}$ }
For the sake of simplicity, we may assume that $q_1 = t_1$.
Then, $(q_1, p_1)$ lies on one of the two exceptional curves $E_1^{\pm}$ at $(t_1, \hat{\nu}_1^{\pm})$.
Suppose that $(q_1, p_1)$ is on $E_1^+$.
We consider the blowing-up of $\widetilde{\mathbb{F}_3}$ at the two points $\{ (q_1, p_1), (q_1, -p_1) \}$.
We denote by $\widetilde{E}_1^+$ the strict transform of $E_1^+$.

In this situation, set
$$ \cK'_{5, q_1} := \Bl_{\{(q_1, p_1), (q_1, -p_1) \}} \widetilde{\mathbb{F}_3} \setminus (\widetilde{s}_{\infty} \cup \widetilde{F}_1 \cup \cdots \cup \widetilde{F}_5 \cup \widetilde{E}_1^+). $$

We will show that the similar result as Corollary $\ref{coh of O}$.
Instead of considering $\cK'_{5, q_1}$, we will consider the following surface:

$$ \cL := \widetilde{\mathbb{F}_3} \setminus (\widetilde{s}_{\infty} \cup \widetilde{F}_1 \cup \cdots \cup \widetilde{F}_5 \cup E_1^+).$$ 

\begin{Prop}\label{coh of L}
\textit{
$$ H^i(\cL, \cO_{\cL}) =   \begin{cases}
                                                                \mathbb{C}, & i = 0, \\
                                                                 0, & i > 0.
                                                      \end{cases} $$
}
\end{Prop}

\begin{proof}
In $\widetilde{\mathbb{F}_3}$, we have that $E_1^+$ is a $(-1)$-curve, and hence we contract this curve. 
Then $\widetilde{F}_1$ becomes a $(-1)$-curve, and we also contract this curve.
As a result, we have the blowing-ups of $\mathbb{F}_2$ at $8$ points, and we have to compute the cohomology of the surface
$$ \cL' := \Bl_{\{8pts\}} \mathbb{F}_2 \setminus (\widetilde{s}_{\infty} \cup \widetilde{F}_2 \cup \cdots \cup \widetilde{F}_5).$$
This is the same situation as \cite[Theorem 2 (iii)]{AL}, and the statement is proved.
\end{proof}

The difference between $\cK'_{5, q_1}$ and $\cL$ is that, adding the points $\{ (q_1, p_1), (q_1, -p_1) \}$, blowing-up these points, and removing the corresponding points.
These operations do not change the cohomology $H^i(\cO)$.

\section{Proof of Theorem \ref{coh of M}}\label{section proof}

By the same argument of Proposition $\ref{second coh of K}$, we have
\begin{Prop}\label{second coh of K5}
\textit{
Let $\cF$ be any quasi-coherent sheaf on $\cK'_{5}$.
Then $H^i(\cK'_{5}, \cF) = 0$ for $i \geq 2$.
}
\end{Prop}

Since $D^{red} := \widetilde{s}_{\infty} + \widetilde{F}_1 + \cdots + \widetilde{F}_5 \subset \widetilde{\mathbb{F}_3}$ is contractible, we have the following lemma.
\begin{Lem}\label{coh of K5}
 \textit{$H^i(\cK'_5, \cO_{\cK'_5}) = \begin{cases}
                                                                                                   \mathbb{C}, & i = 0, \\
                                                                                                    H^2_m(A), & i = 1, \\
                                                                                                                 0, & i \geq 2, 
                                                                                                \end{cases} $\\
where $(A, \mathfrak{m})$ is a local ring such that $\dim(A_{\mathfrak{m}}) = 2$.}
\end{Lem}
\begin{proof}
Let $\pi : \widetilde{\mathbb{F}_3} \ra S$ be a map onto a rational surface $S$ which contracts the divisor $D^{red} \subset \widetilde{\mathbb{F}_3}$ to the rational singular point $\{p \} \subset S$.
Set $U := S \setminus \{ p \}$. Then we have the long exact sequence
\begin{equation*}
 \begin{split}
   0 &\ra H_p^0(S, \cO_S) \ra H^0(S, \cO_S) \ra H^0(U, \cO_U)\\
      &\ra H_p^1(S, \cO_S) \ra H^1(S, \cO_S) \ra H^1(U, \cO_U)\\
      &\ra H_p^2(S, \cO_S) \ra H^2(S, \cO_S) \ra H^2(U, \cO_U)\\
      &\ra H_p^3(S, \cO_S) \ra 0.
 \end{split}
\end{equation*}    
By excision isomorphism, we have $H_p^i(S, \cO_S) = H_p^i(V, \cO_V)$, where $V = \Spec(A)$ and $\{p\}$ corresponds to the maximal ideal $\mathfrak{m}$ of $A$.
Since $V$ is affine, this cohomology is equal to $H_{\mathfrak{m}}^i(A)$.
Now it is straightforward to see that $\dim(A_{\mathfrak{m}}) = \depth_{\mathfrak{m}}(A) = 2$.
Therefore we have $H_{\mathfrak{m}}^i(A) = 0$ for $i \neq 2$, and $H^1(U, \cO_U) \simeq H_{\mathfrak{m}}^2(A) \neq 0$
(see, for example, \cite{H} p.217 exercise 3.4(b)).
\end{proof}

\begin{proof}[Proof of Theorem \ref{coh of M}]
We may assume that $d = -1$.
Set $\widehat{M(-1)}_Z := M(-1)^0 \cup \widehat{Z}$.
By Proposition \ref{injective M0}, we have injective maps $\iota : M(-1)^0 \hookrightarrow \Hilb^2(\cK'_5)$ and $\hat{\iota} : \widehat{M(-1)}_Z \hookrightarrow \widetilde{\Hilb}^2(\cK'_5)$.
We define the blowing-up parameter $\lambda_-$ by $p_1 + p_2= \lambda_- (q_1 - q_2)$.

Set $T := \widetilde{\Hilb}^2(\cK'_5) \setminus \widehat{M(-1)}_Z$. For a vector bundle $\cF$ on $\widetilde{\Hilb}^2(\cK'_5)$,
\begin{equation*}
 \begin{split}
H^i(\widehat{M(-1)}_Z, \cF|_{\widehat{M(-1)}_Z}) &= H^i(\widetilde{\Hilb}^2(\cK'_5), \hat{\iota}_*\hat{\iota}^* \cF) \\ 
                     &= \varinjlim H^i(\widetilde{\Hilb}^2(\cK'_5), \cF(kT)).\\
\end{split}
\end{equation*}
To compute $H^i(\widetilde{\Hilb}^2(\cK'_5), \cF(kT))$, consider $H^i((\cK'_5 \times \cK'_5)^{\approx}, \cF(kT'))$, where 
$T'$ is defined by $ ( \lambda_- = \infty )$.
We can define a map
\begin{equation*}
\begin{aligned}
  f  :  (\cK'_5 \times \cK'_5)^{\approx} \setminus T' &\lra  \cK'_5 \\
 (q_1, p_1, q_2, p_2) &\longmapsto (q_1,p_1),
\end{aligned}
\end{equation*}
and the fiber is $f^{-1}(\{(q_1, p_1)\}) \simeq \cK'_{5, q_1}$.
By Leray's spectral sequence, we have
\begin{equation*}
H^i((\cK'_5 \times \cK'_5)^{\approx} \setminus T', \cF) \simeq \bigoplus_{p + q = i} H^p(\cK'_5, R^qf_* \cF). 
\end{equation*}
Using the base change theorem, we have
$ (R^qf_* \cF)_{(q_1, p_1)} \simeq H^q(\cK'_{5, q_1}, \cF_{(q_1, p_1)})$.
Hence, Theorem \ref{coh of M} follows from Corollary \ref{coh of O}, Lemma \ref{second coh of K5} and Lemma \ref{coh of K5} as follows:
we have
$$H^i((\cK'_5 \times \cK'_5)^{\approx} \setminus T', \cO)  =   \begin{cases} 
                                                               \mathbb{C}, & i = 0, \\
                                                                H^2_{\mathfrak{m}}(A) \oplus H^0(\cK'_5, R^1f_*\cO), & i = 1, \\
                                                                H^1(\cK'_5, R^1f_*\cO), &i=2, \\
                                                                 0, & i > 2\rlap{.}
                                                      \end{cases} $$
Moreover, the action of $\mathfrak{S}_2$ on $H^i((\cK'_5 \times \cK'_5)^{\approx} \setminus T', \cO)$ is nontrivial.
Therefore,
$$H^i(\widehat{M(-1)}_Z, \cO_{\widehat{M(-1)}_Z}) =  \begin{cases}
                                                                                     \mathbb{C}, & i = 0, \\
                                                                                     0, & i > 0\rlap{.}
                                                                              \end{cases}$$
Since $\codim_{\Hilb^2(\cK'_5)}(\widetilde{Z}) = 2$, and $M(-1)^1 = M(-1) \setminus M(-1)^0 \simeq \mathbb{A}^2$ (see \cite{Obl}), we have 
\begin{equation*}
 \begin{split}
 H^i(\widehat{M(-1)}_Z, \cO_{\widehat{M(-1)}_Z})  &= H^i(M(-1)^0 \cup \widetilde{Z}, \cO)\\
                                                                      &= H^i(M(-1)^0, \cO_{M(-1)^0})\\
                                                                      &= H^i(M(-1), \cO_M(-1)) .
  \end{split}
\end{equation*}
\end{proof}

\subsection*{Acknowledgements}
I am very grateful to Professor Masa-Hiko Saito for his constant attention to this work and for his warm encouragement.
I would also like to thank Doctor Arata Komyo for his numerous stimulating discussions and Professor Frank Loray for his hospitality at Universit\'{e} de Rennes 1.


\Address


\begin{thebibliography}{99}
\bibitem{A}
D.\ Arinkin,
{\it Orthogonality of natural sheaves on moduli stacks of SL(2)-bundles with connections on $\mathbb{P}^1$ minus 4 points.},
Selecta Math., New Series 7 (2001), 213-239. 



\bibitem{AL}
D.\ Arinkin, S.\ Lysenko, 
{\it On the moduli of $SL(2)$-bundles with connections on $\bP^1\setminus \{ x_1,\ldots,x_4 \}$}, 
Internat. Math. Res. Notices (1997), no. {\bf 19}, 983--999.

\bibitem{B}
P.\ Boalch,
{\it Geometry of moduli spaces of meromorphic connections on curves, Stokes data, wild nonabelian Hodge theory, hyperkaehler manifolds, isomonodromic deformations, Painleve equations, and relations to Lie theory.},
2012, HAL Id: tel-00768643, https://tel.archives-ouvertes.fr/tel-00768643




\bibitem{BHPV}
W.\ Barth, K.\ Hulek, C.\ Peters, A.\ Van de Ven,
{\it Compact complex surfaces},
Ergebnisse der Mathematik und ihrer Grenzgebiete. 3. Folge. A Series of Modern Surveys in Mathematics [Results in Mathematics and Related Areas. 3rd Series. A Series of Modern Surveys in Mathematics], {\bf 4}, Springer-Verlag, Berlin, 2004.


\bibitem{H}
R.\ Hartshorne,
{\it Algebraic geometry},
Graduate Texts in Mathematics, No. 52, Springer-Verlag, New York-Heidelberg, 1977.



\bibitem{IIS}
M.\ Inaba, K.\ Iwasaki, M.-H.\ Saito, 
{\it Moduli of stable parabolic connections, Riemann-
Hilbert correspondence and geometry of Painlev\'e equation of type VI. I },
Publ. Res. Inst. Math. Sci. (2006), no. {\bf 4}, 987-1089.

\bibitem{LS}
F.\ Loray, M.-H.\ Saito,
{\it Lagrangian fibrations in duality on moduli spaces of rank 2 logarithmic connections over the projective line.} 
Internat. Math. Res. Notices (2015), no. {\bf 4}, 995--1043.


\bibitem{Naka}
H.\ Nakajima,
{\it Lectures on Hilbert schemes of points on surfaces}. University Lecture Series, {\bf18}. American Mathematical Society, Providence, RI, 1999.


\bibitem{KS}
A.\ Komyo, M.-H.\ Saito,
{\it Explicit description of jumping phenomena on moduli spaces of parabolic connections and Hilbert schemes of points on surfaces}, 
Kyoto J. Math. 59 (2019), no. 3, 515–-552.


\bibitem{Obl}
S.\ Oblezin, 
{\it Isomonodromic deformations of $\mathfrak{sl}(2)$ Fuchsian systems on the Riemann sphere.}
Mosc. Math. J. 5 (2005), no. {\bf 2}, 415--441, 494--495.

\bibitem{S}
C.\ Simpson,
{\it An explicit view of the Hitchin fibration on the Betti side for P1 minus 5 points.}
Geometry and physics. Volume 2. Oxford: Oxford University Press (2018), 705--724.


\end{thebibliography}
\end{document}